\documentclass{amsart}
\usepackage{amsmath}
\usepackage{amsfonts}
\usepackage{graphics}
\usepackage{epsfig}
\usepackage{amssymb}
\usepackage{amscd}
\usepackage[all]{xy}
\usepackage{latexsym}
\usepackage{graphicx}
\usepackage{multirow}
\usepackage{longtable}
\usepackage{geometry}

\theoremstyle{remark}

\newtheorem{thm}{\bf Theorem}[section]

\numberwithin{equation}{section} \numberwithin{figure}{section}

\renewcommand*{\to}{\rightarrow}
\renewcommand*{\bar}[1]{\overline{#1}}

\newcommand{\codim}{\operatorname{codim}}

\newcommand{\ind}{\operatorname{ind}}
\newcommand{\Gr}{\operatorname{Gr}}
\newcommand{\End}{\operatorname{End}}
\newcommand{\coker}{\operatorname{coker}}

\newcommand{\ev}{\operatorname{ev}}
\newcommand{\mb}[1]{\mathbb{#1}} % Field

\newcommand{\hs}{\mathcal{H}}% Hilbert space
\newcommand{\mc}[1]{\mathcal{#1}}

\newcommand{\Gor}{\operatorname{Gor}}

\begin{document}
\large \setcounter{section}{0}

\title{Weierstrass cycles and tautological rings in various moduli spaces of algebraic curves}

\author{Jia-Ming (Frank) Liou}
\address{Department of Mathematics, National Cheng Kung University, Taiwan\\fjmliou@mail.ncku.edu.tw}

\author{Albert Schwarz}
\address{Department of Mathematics\\
University of California\\
Davis, CA 95616, USA\\ schwarz@math.ucdavis.edu}

\author{Renjun Xu}
\address{Department of Physics\\
University of California\\
Davis, CA 95616, USA\\rxu@ucdavis.edu}

To Professor Andrzej Granas

\begin{abstract}\large
We analyze  Weierstrass cycles and tautological rings in moduli spaces of smooth algebraic curves
and in moduli spaces of integral algebraic curves with embedded disks with special attention to moduli spaces of  curves having genus $\leq 6.$ In particular, we  show that our general formula gives a good estimate for the dimension of Weierstrass cycles for low genera.
\end{abstract}

\maketitle

\allowdisplaybreaks

%\begin{abstract}
%\end{abstract}

%\maketitle

\allowdisplaybreaks

\section{Introduction}

A numerical semigroup $H$ is a subsemigroup of  $\mb N_{0}$ (of the  semigroup of nonnegative integers)  such that the greatest common divisor of the elements of $H$ is $1$ and $\mb N_{0}\setminus H$ is a finite set. Elements of $\mb N_{0}\setminus H$ are called gaps and the number $g(H)$ of gaps is called the genus of $H.$ Let $C$ be a smooth complex curve of genus $g$ and $p\in C.$ Let $H_{p}$ be the set consisting of $0$ and such integers $n$  that there exists  a function on $C$ holomorphic everywhere except $p$ and having pole of order $n$ at $p.$ Then $H_{p}$ is a numerical semigroup of genus $g$; it is called the Weierstrass semigroup of $p.$ If $H$ is a numerical semigroup of genus $g$ such that there exists a smooth curve $C$ with a point $p$ having Weierstrass semigroup $H,$ we sat that $H$ is  a Weierstrass semigroup.

The moduli space of complex smooth pointed curves\footnote{All the curves in this paper are assumed to be irreducible projective curves.} of genus $g$ is denoted by $\mc M_{g,1}$ and has dimension\footnote{The dimensions and codimensions of spaces considered in this paper are over $\mb C.$} $3g-2$ for $g>0.$ If $H$ is a numerical semigroup of genus $g,$ we denote $\mc M_{H}$ the subspace of $\mc M_{g,1}$ consisting of points $(C,p)$ such that $p$ has Weierstrass semigroup $H.$ (This subspace can be regarded as locally closed subscheme.) One can say that   $\mc M_{H}$  is the moduli space of smooth pointed curves with prescribed Weierstrass semigroup $H$. Some interesting questions arise when studying $\mc M_{H}.$ When is the space $\mc M_{H}$ nonempty? If it is nonempty, what is its dimension (or codimension in $\mc M_{g,1}$)? In \cite{EH}, Eisenbud and Harris obtained an lower bound of $\dim\mc M_{H}$ for any irreducible component of $\mc M_{H}$ (see (\ref{EH}) below); in \cite{De}, Deligne obtained the estimate (\ref{del}) which gives the upper bound of $\dim\mc M_{H}.$ It was shown  in \cite{JK} that any numerical semigroup of genus $\leq 7$ is a Weierstrass semigroup. In \cite{NM} and \cite{Nak}, Nakano and Mori studied the dimension of $\mc M_{H}$ for $g\leq 6.$

The closure of $\mc M_{H},$ denoted by $W_{H},$ is called a Weierstrass cycle of semigroup $H$ (or \textbf{W}-cycle for short). In \cite{LS3},  using Krichever map, we have found a general estimate from below for the dimension of Weierstrass cycles. In the case when this estimate is precise we  expressed the cohomology classes  dual to the cycles $W_{H}$ in terms of $\lambda$ and $\psi$-classes.  (This result is based  on the study of cohomological properties of Krichever map in \cite {LS}.) In present paper we apply  these results to low genera. In section 2, we discuss the estimate for the dimension of Weierstrass cycles in more detail; we  show that it is precise for $g\leq 5$  and for $g=6$ the error is $\leq 1.$  In section 4, we list the cohomology classes of $W_{H}$ for $g\leq 6.$ 

The definition of Weierstrass point  can be applied also to singular curves if $p$ is a smooth point.
Moreover, all our constructions can be generalized to the case of integral curves (= irreducible reduced projective curves). In particular, we can define Weierstrass cycles in moduli spaces of integral curves. We will show that for small genera the dimensions of these cycles coincide with dimensions of Weierstrass cycles in the moduli spaces of smooth curves.

 In \cite {LS3} we have used the identification of the moduli space $\widehat{\mathcal{CM}}_{g}$ of complex integral curves of genus $g$ with embedded disks with Krichever locus of the Segal-Wilson version the Sato Grassmannian. This moduli space can be represented also by a closed subscheme of the Sato Grassmaniann $\Gr(\hs),$ defined in the framework of algebraic geometry \cite{Mul}, \cite{MP}.
The circle $T=U(1)$ acts on both of $\widehat{\mathcal{CM}}_{g}$ and $\Gr(\hs);$ the Krichever embedding $k$ of $\widehat{\mathcal{CM}}_{g}$ into $\Gr(\hs)$ induces a homomorphism $k^*$ of equivariant cohomology. The tautological ring of $\widehat{\mathcal{CM}}_{g}$ can be defined as the image of $k^*$ or as the smallest graded complex subalgebra $R=\bigoplus_{i\geq 0}R_{i}$ of $H_{T}^{*}(\widehat{\mathcal{CM}}_{g})$ generated by the equivariant $\lambda$ and $\psi$-classes (see Section 3 or \cite {LS3} for the definition of these classes). In section 3, we give explicit formulas for relations in tautological ring of $\widehat{\mathcal{CM}}_g$ for $g\leq 6.$ We do not claim that all relations in this ring follow from our relations.  Our relations give only an estimate from above for the size of this ring. We obtain an estimate from below considering the restriction of cohomology classes to fixed points of $T$-action. (These fixed points correspond to monomial curves.)  We can define  the set $\widehat {\mathcal{CM}_H}$  in  $\widehat{\mathcal{CM}_{g}}$  fixing the Weierstrass semigroup $H$ at the marked point and its closure, the Weierstrass cycle $\widehat{\mathcal{C} {W_H} }$.  In the subspace of smooth curves with disks $\widehat{\mathcal{M}_g}\subset\widehat{\mathcal{CM}}_g$ we have corresponding objects  $\widehat {\mc{M}_H}$ and $\widehat {W_H}.$ The equivariant cohomology classes $[W_H]$ corresponding to $T$-invariant cycles  $\widehat {W_H}$ can   be identified with cohomology classes in  $\mc M_{g,1}$ corresponding to Weierstrass cycles $W_H$
(see \cite {LS3}). We study  classes $[W_H]$  in equivariant cohomology; for smooth curves this gives us information about classes of Weierstrass cycles in ordinary cohomology.

To facilitate the reading  we review some of results of \cite {LS3}. The archive version of the present paper  arXiv:1308.6374  contains some details that are  skipped in the journal version.

In \cite{Sto}, St\"{o}hr constructed the moduli spaces $\mc M_{H}^{\Gor}$ of Gorenstein curves with symmetric Weierstrass semigroup $H$; in \cite{CS}, Contiero and St\"{o}hr studied the dimension of $\mc M_{H}^{\Gor}.$ In the cases when their calculations and our calculations can be compared they agree. 

%In section 5  we will construct the moduli space $\mc{CM}_H$  of integral curves with a prescribed numerical semigroup $H$ of genus $g$ at the marked point; it is closely related to the  Weierstrass cycle $\widehat{\mc {C} {W_H}}. $

\section{Estimates}
The Krichever map allows us to estimate from below of the dimension of Weierstrass cycles (simply \textbf{W}-cycles).     We compare our estimate with the results obtained from the moduli space of curves of low genera in \cite {NM}, \cite {Nak} and \cite {CS}. For genus $\leq 5,$ we see that our estimate is precise; together with the irreducibility of all \textbf{W}-cycles (except one) proven in \cite {Nak}, this permits us to calculate the corresponding cohomology classes (up to a constant factor). For genus $6$, Nakano calculated the dimension of \textbf{W}-cycles and proved irreducibility of $16$ \textbf{W}-cycles out of $23$. Our estimate together with Deligne estimate from above gives the exact answer for $6$ cycles that cannot be analyzed by methods of \cite {Nak}. There is one \textbf{W}-cycle whose dimension cannot be found exactly neither by methods of \cite {Nak}, nor by our methods. We were able to prove that the codimension is either $2$ or $3.$ In $22$ cases when the dimension is known, our estimate is precise in $16$ cases; in the remaining cases, the difference between the exact dimension of \textbf{W}-cycles and our estimate is equal to $1$.

The \textbf{W}-cycle $W_{H}$ is related to the Schubert cycle $\bar{\Sigma}_{S}$ on the Grassmannian via the Krichever map.  As we mentioned already we can work either with Segal-Wilson definition  of Grassmannian or in algebraic approach; let us review some  basic definitions in algebraic approach . (For details, see \cite{MP} and \cite{Mul}. ) The field of Laurent series $\mb C((z))$ is a complete complex topological vector space (field) with respect to the filtration $\{z^{n}\mb C[[z]]:n\in\mb Z\}$. Let $\hs$ be the closed linear subspace of $\mb C((z))$ generated by $\{z^{i}:i\neq -1\}$ together with a direct sum decomposition $\hs=\hs_{-}\oplus\hs_{+},$ where $\hs_{+}=\mb C[[z]]$ and $\hs_{-}=z^{-2}\mb C[z^{-1}].$ The Sato Grassmannian $\Gr(\hs)$  in this approach is a connected reduced separated scheme over $\mb C$ whose points are closed subspaces $W$ of $\hs$ such that the canonical projection $\pi_{W}:W\mapsto \hs_{-}$ is Fredholm whose index $\ind\pi_{W}=\dim\ker\pi_{W}-\dim\coker\pi_{W}$ is $g.$

The Schubert cells $\Sigma_{Z}$ of $\Gr(\hs)$ are labeled by sequences of virtual cardinality $g.$ Here a sequence of virtual cardinality $g$ is a decreasing sequence of integers $Z$ with $z_{i}\neq -1$ for all $i$ such that both of $Z\cap \mb Z_{\geq 0}$ and $\mb Z_{<0}-Z$ are finite sets, and $\#(Z\cap \mb Z_{\geq 0})-\#(\mb Z_{<0}-Z)=g$. One can check that if $Z=(z_{i})$ is a sequence of virtual cardinality $g,$ then $z_{i}=-i+g-1$ for $i\gg 0.$ In this case, the codimension of the Schubert cell $\Sigma_{Z}$ in $\Gr(\hs)$ is given by the formula
\begin{equation*}
w(Z)=\sum_{i=1}^{i_{0}}(z_{i}+i-g)+\sum_{i=i_{0}}^{\infty}(z_{i}+i-g+1),
\end{equation*}
where $i_{0}$ is the index such that $z_{i_{0}}\geq 0$ and $z_{i_{0}+1}<0.$ 

The Krichever map
\begin{equation*}
k:\widehat{\mathcal{CM}}_{g}\to\Gr(\hs)
\end{equation*}
is a closed immersion defined by $(C,p,z)\mapsto z(H^{0}(C- p,\omega_{C})),$ where $\omega_{C}$ is the dualizing sheaf of $C$ on $\omega_{C}.$ Moreover, the moduli space $\widehat{\mc M}_{g}$ is contained in $\widehat{\mathcal{CM}}_{g}$ as an open subscheme.

Let $F:\widehat{\mathcal{M}}_{g}\to\mc M_{g,1}$ be the forgetul map $(C,p,z)\mapsto (C,p).$ We see that
\begin{equation}\label{codim}
\codim(\widehat{W}_{H};\widehat{\mc M}_{g})=\codim(W_{H};\mc M_{g,1}).
\end{equation}
 Here we use the notation $\widehat{W}_{H}=F^{-1}W_{H}$ for any Weierstrass semigroup $H.$
Later we will denote  $\codim(W_{H};\mc M_{g,1})$ by $\codim W_{H}$ for simplicity.

Let $\varsigma:\mb Z\to\mb Z$ be the translation operator $n\mapsto n+1$. Let us denote by $S$ the   set $\varsigma^{-1}(\mb Z-H)$  where $H$ is  a Weierstrass semigroup .  The set $S$ can be considered as a decreasing sequence\footnote{Let $s_{1}$ be the largest integer in $S$ and define $s_{j}$ recursively by requiring that $s_{j}$ is the largest integer in $S\setminus\{s_{1},\cdots,s_{j-1}\}.$ } of integers that is called the Weierstrass sequence corresponding to $H.$ It follows from the Riemann-Roch theorem that every Weierstrass sequence $S$ is a sequence of virtual cardinality $g$ with $s_{g}=0.$ One can check that the \textbf{W}-cycle $\widehat{W}_{H}$ is the preimage of $\bar{\Sigma}_{S}$ under the Krichever map, i.e. $\widehat{W}_{H}=k^{-1}(\bar{\Sigma}_{S}).$ In other words the Weierstrass cycles are intersections of Schubert cells with the Krichever locus. Notice that $\widehat{W}_{H} $ is a closure of $\widehat{M}_{H}=k^{-1}({\Sigma}_{S})$ that can be regarded as a moduli space of smooth curves with embedded disks and with prescribed Weierstrass semigroup $H$ at the center of the disk.

Let $\Sigma_{Z}$ be a Schubert cell in $\Gr(\hs).$ The W-cycle $\widehat{W}_{Z}$ is a union of W cycles $\widehat{W}_S$ with $S\geq Z$. It is obvious that this cycle (an intersection  of $\bar{\Sigma}_{S}$ with the Krichever locus) has codimension $\leq w(Z)$. By (\ref{codim}) and taking into account that the codimension of a union is the minimal codimension of components, we obtain the following statement:
\begin {thm}
If $S$ runs over Weierstrass sequences obeying $S\geq Z,$ then
\begin{equation}\label{est1}
\min \codim W_{S}\leq w(Z).
\end{equation}
\end {thm}
The codimension of the Schubert cycle $\bar{\Sigma}_{S}$ labeled by a Weierstrass sequence $S$ is $\sum_{i=1}^{g}(s_{i}+i-g).$ We see that (\ref{est1}) is a stronger statement than Eisenbud-Harris (EH) estimate
\begin{equation}\label{EH}
\codim W_{S}\leq w(S)
\end{equation}
(see \cite {EH}). Note that the relation between the Weierstrass sequence $S$ and the Weierstrass gap sequence $\Gamma=\{1=n_{1}<\cdots<n_{g}\}$ of a Weierstrass semigroup $H$ is given by $s_{i}=n_{g-i+1}-1$ for $1\leq i\leq g.$ Since a Weierstrass sequence $S$ is determined by its first $g$-values and the set of its first $g$-values is in one-to-one correspondence with the Weierstrass gap sequence $\Gamma$, we denote $W_{S}$ by $W_{\Gamma}.$

If $S,Z$ are two decreasing sequences of integers of virtual cardinality $g,$ we say that $S\leq Z$ iff $s_{i}\leq z_{i}$ for all $i.$ Hence if $S$ is a Weierstrass sequence and $Z$ a sequence of virtual cardinality $g$ obeying $Z\leq S$, then $z_{i}=g-i-1$ for $i\geq g+1.$ From now on, we only consider sequences $Z$ obeying $z_{i}=g-i-1$ for $i\geq g+1.$ In this case, the finite increasing sequence $\{z_{g}+1\leq \cdots\leq z_{1}+1\}$ is called the ``gap sequence'' for $Z.$

We are using  the description of Weierstrass gap sequences for $g\leq 6$ given in \cite{Nak} and \cite{NM}.

For $g\leq 3,$ the only case when the EH estimate is not exact is the case of Weierstrass sequence $S$ with the gap sequence $\{1,3,5\}$. We take $Z$ with the gap sequence $\{1,3,4\}.$ The Weierstrass sequence obeying $S\geq Z$ is unique. We obtain $\codim W_{1,3,5}\leq 2.$

%Checked

For $g=4,$ the EH estimate is not exact only in the cases of gap sequences $\{1,3,5,7\}$ and $\{1,2,4,7\}$; we use $Z$ with gap sequences $\{1,3,4,5\}$ and $\{1,2,4,6\}$ and obtain that $\codim W_{1,3,5,7}\leq 3$ and that $\codim W_{1,2,4,7}\leq 3.$

%Checked

For $g=5$ the EH estimate is not exact only in the cases of gap sequences $\{1,3,5,7,9\}$ and $\{1,2,3,5,9\}$; we use $Z$ with gap sequences $\{1,3,4,5,6\}$ and $\{1,2,3,5,8\}$; we see that $\codim W_{1,3,5,7,9}\leq 4$ and $\codim W_{1,2,3,5,9}\leq 4$.

For $g=6,$ the EH estimate is not exact for the following gap sequences:$\{1,2,3,5,7,11\}$, $\{1,2,3,6,7,11\}$, $\{1,2,4,5,7,8\}$, $\{1,2,4,5,7,10\}$, $\{1,2,4,5,8,11\}$, $\{1,3,5,7,9,11\}$. Together with Deligne's estimate,
\begin{equation}\label{del}
\dim \mc M_{H}\leq 2g-[\End H:H]-2,
\end{equation}
we determine the dimensions of $\mc M_{H}$ which were not found in \cite{Nak}. The results are listed in the following table.
\begin{center}
\begin{tabular}{|l|c|c|l|c|}
  \hline
  gap sequence of $S$  & E-H & our estimate & gap sequence of $Z$ & exact codim \\ \hline
      $\{1,2,3,5,7,11\}$ & $8$ & 7 &$\{1,2,3,4,7,11\}$ & $6$\\
      $\{1,2,3,6,7,11\}$ & $9$ & 6 & $\{1,2,3,6,7,8\}$ & $6$\\
      $\{1,2,4,5,7,8\}$  & $6$ & 4 & $\{1,2,4,5,6,7\}$ & $4$\\
      $\{1,2,4,5,7,10\}$ & $8$ & 6 & $\{1,2,4,5,6,9\}$ & $5$\\
      $\{1,2,4,5,8,11\}$ &$10$ & 7 & $\{1,2,3,5,8,9\}$ & $6$\\
      $\{1,3,5,7,9,11\}$ &$15$ & 5 & $\{1,3,4,5,6,7\}$ & $5$\\
  \hline
\end{tabular}
\end{center}
\section{Tautological Ring of $\widehat{\mathcal{CM}}_{g}$}
The circle $T=U(1)$ acts on both of $\widehat{\mathcal{CM}}_{g}$ and $\Gr(\hs).$ Denote $u=c_{1}(\mc O_{\mb P^{\infty}}(1))=c_1(B_T)$ the first Chern class of the tautological line bundle over $\mb P^{\infty}$ (over the classifying space $B_T.$) Then both of $H_{T}^{*}(\Gr(\hs))$ and $H_{T}^{*}(\widehat{\mathcal{CM}}_{g})$ are modules over $\mb C[u].$ The Krichever map $k:\widehat{\mathcal{CM}}_{g}\to \Gr(\hs)$ is equivariant and induces a morphism $k^{*}:H_{T}^{*}(\Gr(\hs))\to H_{T}^{*}(\widehat{\mathcal{CM}}_{g}).$

The class $\psi=-k^{*}u$ in $H_{T}^{*}(\widehat{\mathcal{CM}}_{g})$ is called the equivariant $\psi$-class. The Hodge bundle $\mb E$ over $\widehat{\mathcal{CM}}_{g}$ is a $T$-equivariant vector bundle of rank $g$ whose fiber over $(C,p,z)$ is the vector space $z(H^{0}(C,\omega_{C})).$ The $i$-th equivariant Chern class of $\mb E$ is denoted by $\lambda_{i}$ for $1\leq i\leq g.$ (We will keep the same notation for $\psi$ and $\lambda$-classes when we consider their images in $H_{T}^{*}(\widehat{\mc M}_{g})$).

The tautological ring $R=R(\widehat{\mathcal{CM}}_{g})$ is the $\mb Q$-subalgebra of $H_{T}^{*}(\widehat{\mathcal{CM}}_{g})$ generated by $\lambda_{i},$ $1\leq i\leq g,$ and $\psi.$ We assign $\deg\lambda_{i}=i$ and $\deg \psi=1.$ Then $R$ is a graded $\mb Q$-algebra $\bigoplus_{i\geq 0}R_{i}$. Set $h_{i}(R)=\dim_{\mb Q}R_{i},$ for all $i.$ The Hilbert Poincare series of $R$ is denoted by $P_{R}(t)=\sum_{i}h_{i}(R)t^{i}\in\mb Z[[t]].$

Consider the free polynomial algebra $A=\mb Q[\lambda_{1},\cdots,\lambda_{g},\psi]$ generated by the symbols $\lambda_{1},\cdots,\lambda_{g},$ and $\psi$ with $\deg \lambda_{i}=i$ and $\deg \psi=1.$ Then $R=A/I_{tau}$ for the ideal $I_{tau}$ called the ideal of tautological relations.

Since the monomial curves are in one-to-one correspondence with the $T$-fixed points on $\widehat{\mathcal{CM}}_{g},$ we obtain a ring homomorphism $\ev$ from $R$ to $\bigoplus_{S}\mb C[\psi]$. Here $S$ runs over all Weierstrass sequence. The ideal $I_{tau}$ is contained in the kernel $I_{\ev}$ of $\ev.$ This gives us a surjective ring homomorphism $R\to A/I_{ev}.$ On the other hand, if the intersection of $\bar{\Sigma}_{\mu}$ and $k(\widehat{\mathcal{CM}}_{g})$ is empty, then $k^{*}\Omega_{\mu}^{T}=0$ in $R.$ We denote by $I$ the ideal of $A$ generated by $k^{*}\Omega_{\mu}^{T}.$ Then $I\subset I_{tau}.$ Then we obtain another surjective ring homomorphism $A/I\to R.$ As a consequence, we have the following estimates
\begin{equation*}
h_{i}(A/I_{\ev})\leq h_{i}(R)\leq h_{i}(A/I).
\end{equation*}
The Hilbert Poincare series of $A/I_{ev}$ and $A/I$ are given by the following table for $2\leq g\leq 6$:

\begin{equation*}
\begin{tabular}{|l|l|l|}
  \hline
  $g$ & $P_{A/I_{ev}}$& $P_{A/I}$\\ \hline
  2 & $1+t$ & $(1+t+2t^{2}-2t^{4}-t^{5})(1-t)^{-1}$  \\ \hline
  3 & $1+2t+2t^{2}+t^{3}$ & $1+2t+4t^{2}+7t^{3}+9t^{4}+9t^{5}+6t^{6}+t^{7}$ \\ \hline
  4 & $1+2t+4t^{2}+3t^{3}+t^{4}$ & $1+2t+4t^{2}+7t^{3}+12t^{4}+16t^{5}+20t^{6}$ \\
   & &  $+22t^{7}+21t^{8}+15t^{9}+9t^{10}+2t^{11}$  \\ \hline
  5 & $1+2t+4t^{2}+7t^{3}+2t^{4}+t^{5}$ & $1+2t+4t^{2}+7t^{3}+12t^{4}+19t^{5}+27t^{6}$ \\
  & & $+35t^{7}+43t^{8}+51t^{9}+54t^{10}+54t^{11}$\\
  & & $+49t^{12}+41t^{13}+27t^{14}+12t^{15}+2t^{16}$  \\ \hline
  6 & $1+2t+4t^{2}+7t^{3}+11t^{4}+6t^{5}+3t^{6}$ & $1+2t+4t^{2}+7t^{3}+12t^{4}+19t^{5}+30t^{6}$ \\
  & & $+42t^{7}+57t^{8}+73t^{9}+92t^{10}+110t^{11}$\\
  && $+127t^{12}+138t^{13}+149t^{14}+151t^{15}+144t^{16}$\\
  && $+129t^{17}+106t^{18}+75t^{19}+41t^{20}+15t^{21}+2t^{22}$
  \\ \hline
\end{tabular}
\end{equation*}

\section{Cohomology Classes of \textbf{W}-cycles}

We are interested in equivariant cohomology classes $[W_{S}]$ corresponding to $T$-invariant  \textbf{W}-cycles $W_{S}.$
If   \textbf{W}-cycle $W_{S}$ is irreducible and our estimate of its dimension is precise  Theorem 2.4 of \cite {LS3} allows us to calculate $[W_{S}]$ up to a constant factor. More rigorously we can formulate this statement in the following way. We assume that we have an equality in (\ref{est1}). In other words, we assume that $\codim W_{S}=w(Z)$ and  for all Weierstrass sequences $\tilde S$ such that $\tilde S\geq Z$ we have $\codim W_{\tilde S}>w(Z)$ . (Notice, that we have made an additional assumption that the sequence obeying $\codim W_{S}=w(Z)$  is unique.) We can say  that  in this situation the intersection of  Schubert cycle  $\bar{\Sigma}_{Z}$ with the Krichever locus is equal to the \textbf{W}-cycle $W_{S}$ and the intersection has expected codimension.  This  allows  us to use Theorem 2.4 of \cite {LS3} to obtain the expression for $[W_{S}]$ in terms of $\psi$ and $\lambda _i.$ We get concrete expression for $[W_{S}]$ in the case when the genus $\leq 5$  calculating the necessary factorial Schur functions (see the Appendix to the archive version of present paper for more detail).

When $g=2,3,4$ all the Weierstrass cycles has the expected codimension. When $g=5,$ the Weierstrass cycles $W_{(2,1,1)},$ $W_{1,2,4,5,8}=W_{(3,1,1)}$ and $W_{1,2,3,5,9}=W_{(4,1)}$ do not have the expected codimension. Moreover, $W_{(2,1,1)}=W_{(1,1,1)}.$ We list results in the case when $W_{\mu}$ has the expected codimension in the following table. (Let us emphasize  that we are calculating the cohomology classes only up to a constant factor.)

\begin{center}
\begin{table}[h]
\caption[${[W_{\mu}]}$]{${[W_{\mu}]}$ for genus $2-5$}
\begin{tabular}{|c|l|c|c|l|}
  \hline
  genus &$W_{\mu}$ & codim & $|\mu|$ & ${[W_{\mu}]}$ \\ \hline
  2& $W_{(1)}$ & 1 & 1 & $3\psi-\lambda_{1}.$\\ \hline
   &$W_{(1)}$& 1 & 1 & $6\psi-\lambda_{1}$ \\
  3&$W_{(1,1)}$ & 2 & 2 & $7\psi^{2}-3\psi\lambda_{1}+\lambda_{1}^{2}-\lambda_{2}$ \\
  &$W_{(2)}$ & 2 & 2 & $35\psi^{2}-10\psi\lambda_{1}+\lambda_{2}$\\
  \hline
  &$W_{(1)}$ &1 &1 & $10\psi-\lambda_{1}$\\
  &$W_{(1,1)}$&2&2 & $25\psi^{2}-6\psi\lambda_{1}+\lambda_{1}^{2}-\lambda_{2}$\\
  4& $W_{(2)}$   &2&2 & $85\psi^{2}-15\psi\lambda_{1}+\lambda_{2}$\\
  &  $W_{(1,1,1)}$    &3 &3  &$15\psi^{3}+3\psi\lambda_{1}^{2}-\lambda_{1}^{3}-3\psi\lambda_{2}+\lambda_{1}(-7\psi^{2}+2\lambda_{2})-\lambda_{3}$ \\
  & $W_{(2,1)}$ & 3 & 3 & $285\psi^{3}+15\psi\lambda_{1}^{2}-9\psi\lambda_{2}-\lambda_{1}(90\psi^{2}+\lambda_{2})+\lambda_{3}$\\
  & $W_{(3)}$& 3 & 3 & $735\psi^{3}-175\psi^{2}\lambda_{1}+21\psi\lambda_{2}-\lambda_{3}$\\ \hline
  &$W_{(1)}$ & 1 & 1 & $15\psi-\lambda_{1}$ \\
  &$W_{(1,1)}$ & 2 & 2 & $65\psi^{2}-10\psi\lambda_{1}+\lambda_{1}^{2}-\lambda_{2}$ \\
  &$W_{(2)}$ & 2 & 2 & $175\psi^{2}-21\psi\lambda_{1}+\lambda_{2}$ \\
  5 &$W_{(1,1,1)}$ & 3 & 3 & $90\psi^{3}+6\psi\lambda_{1}^{2}-\lambda_{1}^{3}-6\psi\lambda_{2}+\lambda_{1}(-25\psi^{2}+2\lambda_{2})-\lambda_{3}$ \\
  &$W_{(2,1)}$ & 3 & 3 &$1015\psi^{3}+21\psi\lambda_{1}^{2}-11\psi\lambda_{2}-\lambda_{1}(210\psi^{2}+\lambda_{2})+\lambda_{3}$ \\
  &$W_{(3)}$ & 3 & 3 & $1960\psi^{3}-322\psi^{2}\lambda_{1}+28\psi\lambda_{2}-\lambda_{3}$  \\
  & $W_{(1,1,1,1)}$ & 4 & 4 & $31\psi^{4}-15\psi^{3}\lambda_{1}+7\psi^{2}\lambda_{1}^{2}-3\psi\lambda_{1}^{3}+\lambda_{1}^{4}-7\psi^{2}\lambda_{2}$\\
  & & & &$+6\psi\lambda_{1}\lambda_{2}-3\lambda_{1}^{2}\lambda_{2}+\lambda_{2}^{2}-3\psi\lambda_{3}+2\lambda_{1}\lambda_{3}-\lambda_{4}$ \\
  &$W_{(2,2)}$ & 4 & 4 & $3850\psi^{4}-1050\psi^{3}\lambda_{1}+140\psi^{2}\lambda_{1}-55\psi^{2}\lambda_{2}-15\psi\lambda_{1}\lambda_{2}$\\
  && & & $+\lambda_{2}^{2}+15\psi\lambda_{3}-\lambda_{1}\lambda_{3}$ \\
  &$W_{(3,1)}$ & 4 & 4 &  $12831\psi^{4}-3220\psi^{3}\lambda_{1}+322\psi^{2}\lambda_{1}^{2}-42\psi^{2}\lambda_{2}$\\
  &&&& $-28\psi\lambda_{1}\lambda_{2}+18\psi\lambda_{3}+\lambda_{1}\lambda_{3}-\lambda_{4}$\\
  &$W_{(4)}$ &4 & 4& $22449\psi^{4}-4536\psi^{3}\lambda_{1}+546\psi^{2}\lambda_{2}-36\psi\lambda_{3}+\lambda_{4}$ \\
  \hline
\end{tabular}
\end{table}
\end{center}

If we consider  Weierstrass cycles in $\widehat{\mathcal{CM}}_{g}$ then instead of  Weierstrass sequences we should use sequences corresponding to numerical semigroups. ( Every numerical semigroup corresponds to a Weierstrass point of a singular curve;  a similar statement is wrong for smooth curves. However, for small genera ($g<7$) this statement is correct \cite {JK}.) Notice, however, that in this case the definition of the cohomology class corresponding to a cycle should be modified: we define this cohomology class as a class containing a cocycle with a support on the given cycle. (Of course, this definition does not specify  the cohomology class uniquely; in our examples it is defined up to a constant factor.)

\section{Appendix}

The class $[W_{\mu}]$ of Weierstrass cycles can be expressed in terms of factorial Schur functions (up to a constant factor). The factorial Schur polynomial $^{n}t_{\mu}(z_{1},\cdots,z_{n})$ of partition $\mu$ in variables $\{z_{1},\cdots,z_{n}\}$
is given by the formula:
\begin{equation*}
^{n}t_{\mu}(z_{1},\cdots,z_{n})=\frac{\det[(z_{i}\downharpoonright
\mu_{j}+n-j)]_{i,j=1}^{n}}{\det[(z_{i}\downharpoonright n-j)]_{i,j=1}^{n}},
\end{equation*}
where the symbol $(z\downharpoonright i)$ stands for the $i$-th
falling factorial power of the variable $z$:
\begin{equation*}
(z\downharpoonright i)=\left\{
                         \begin{array}{ll}
                           z(z-1)\cdots(z-i+1), & \hbox{$i=1,2,\cdots$;} \\
                           1, & \hbox{$i=0$.}
                         \end{array}
                       \right.
\end{equation*}

\begin{thm}
\cite{LS3} If the complex codimension of $W_S$ is equal to $|\mu|=\sum \mu_i$ then
\begin{equation}
\label{w}
[W_S]=\mbox{const}\ ^{g}t_{\mu}\left(-\frac{x_{1}}{\psi},\cdots.-\frac{x_{g}}{\psi}\right)(-\psi)^{|\mu|},
\end{equation}
Here $\{x_{1},\cdots,x_{g}\}$ are the Chern roots of the Hodge bundle $\mb E$ and $\mu$ is the partition corresponding to the sequence $S.$
\end{thm}

The factorial Schur function $^{g}t_{\mu}$ is an inhomogeneous symmetric polynomial; we represent it as a sum of homogeneous symmetric polynomials
$$^{g}t_{\mu}\left(z_{1},\cdots,z_{g}\right)=\sum_{i}t_{\mu}^{i}\left(z_{1},\cdots,z_{g}\right),$$
where $t_{\mu}^{i}(z_{1},\cdots,z_{g})\in\mb C[z_{1},\cdots,z_{g}]^{S_{g}}$ is homogeneous polynomial of degree $i$ . Then (\ref{w}) can be rewritten as
$$[W_{S}]=\sum_{i}t_{\mu}^{i}(x_{1},\cdots,x_{g})(-\psi)^{|\mu|-i}.$$
In the table below, we listed the homogeneous component $t_{\mu}^{i}$  for  $2\leq g\leq 5.$ 

\begin{center}
%	\begin{tabular}{|c|l|l|}
\begin{longtable}{|c|l|p{.75\textwidth}|}
		\caption{$t_{\mu}^i \, (i=0,\cdots, |\mu|)$ for  $2\leq g \leq5$}\\
		\hline
		genus &$t_{\mu}^i$ & $ t_{\mu}^i(x_1,\cdots,x_g)$ \\ 
		\hline
		\endfirsthead
		\multicolumn{3}{c}%
		{\tablename\ \thetable\ -- \textit{Continued from previous page}} \\
		\hline
		genus &$t_{\mu}^i $ & $ t_{\mu}^i(x_1,\cdots,x_g)$ \\ 
		\hline
		\endhead
		\hline \multicolumn{3}{r}{\textit{Continued on next page}} \\
		\endfoot
		\hline
		\endlastfoot
		2& $t_{(1)}^0$ & $-3$\\ 
		& $t_{(1)}^1$ & $x_1+x_2$\\
		\hline
		&$t_{(1)}^0$ & $-6$ \\
		&$t_{(1)}^1$ & $x_1+x_2+x_3$ \\
		&$t_{(1,1)}^0$ & $7$ \\
		&$t_{(1,1)}^1$ & $-3 x_1-3 x_2-3 x_3$ \\
		3&$t_{(1,1)}^2$ & $x_1 x_2+x_1 x_3+x_2 x_3$ \\
		&$t_{(2)}^0$ & $35$\\
		&$t_{(2)}^1$ & $-10 x_1-10 x_2-10 x_3$\\
		&$t_{(2)}^2$ & $x_1^2+x_1 x_2+x_2^2+x_1 x_3+x_2 x_3+x_3^2$\\
		\hline
		&$t_{(1)}^0$ & $-10$\\
		&$t_{(1)}^1$ & $x_1+x_2+x_3+x_4$\\
		&$t_{(1,1)}^0$ & $25$\\
		&$t_{(1,1)}^1$ & $-6 x_1-6 x_2-6 x_3-6 x_4$\\
		&$t_{(1,1)}^2$ & $x_1 x_2+x_1 x_3+x_2 x_3+x_1 x_4+x_2 x_4+x_3 x_4$\\
		4& $t_{(2)}^0$ & $85$\\
		& $t_{(2)}^1$ & $-15 x_1-15 x_2-15 x_3-15 x_4$\\
		& $t_{(2)}^2$ & $x_1^2+x_1 x_2+x_2^2+x_1 x_3+x_2 x_3+x_3^2+x_1 x_4+x_2 x_4+x_3 x_4+x_4^2$\\
		&  $t_{(1,1,1)}^0$  &$-15$ \\
		&  $t_{(1,1,1)}^1$  &$7 x_1+7 x_2+7 x_3+7 x_4$ \\
		&  $t_{(1,1,1)}^2$  &$-3 x_1 x_2-3 x_1 x_3-3 x_2 x_3-3 x_1 x_4-3 x_2 x_4-3 x_3 x_4$ \\
		&  $t_{(1,1,1)}^3$  &$x_1 x_2 x_3+x_1 x_2 x_4+x_1 x_3 x_4+x_2 x_3 x_4$ \\
		& $t_{(2,1)}^0$ & $-285$\\
		& $t_{(2,1)}^1$ & $90 x_1+90 x_2+90 x_3+90 x_4$\\
		& $t_{(2,1)}^2$ & $-6 x_1^2-21 x_1 x_2-6 x_2^2-21 x_1 x_3-21 x_2 x_3-6 x_3^2-21 x_1 x_4-21 x_2 x_4-21 x_3 x_4-6 x_4^2$\\
		& $t_{(2,1)}^3$ & $x_1^2 x_2+x_1 x_2^2+x_1^2 x_3+2 x_1 x_2 x_3+x_2^2 x_3+x_1 x_3^2+x_2 x_3^2+x_1^2 x_4+2 x_1 x_2 x_4+x_2^2 x_4+2 x_1 x_3 x_4+2 x_2 x_3 x_4+x_3^2 x_4+x_1 x_4^2+x_2 x_4^2+x_3 x_4^2$\\
		& $t_{(3)}^0$ & $-735$\\ 
		& $t_{(3)}^1$ & $175 x_1+175 x_2+175 x_3+175 x_4$\\ 
		& $t_{(3)}^2$ & $-21 x_1^2-21 x_1 x_2-21 x_2^2-21 x_1 x_3-21 x_2 x_3-21 x_3^2-21 x_1 x_4-21 x_2 x_4-21 x_3 x_4-21 x_4^2$\\ 
		& $t_{(3)}^3$ & $x_1^3+x_1^2 x_2+x_1 x_2^2+x_2^3+x_1^2 x_3+x_1 x_2 x_3+x_2^2 x_3+x_1 x_3^2+x_2 x_3^2+x_3^3+x_1^2 x_4+x_1 x_2 x_4+x_2^2 x_4+x_1 x_3 x_4+x_2 x_3 x_4+x_3^2 x_4+x_1 x_4^2+x_2 x_4^2+x_3 x_4^2+x_4^3$\\ 
		\hline
		&$t_{(1)}^0$ & $-15$ \\
		&$t_{(1)}^1$ & $x_1+x_2+x_3+x_4+x_5$ \\
		&$t_{(1,1)}^0$ & $65$ \\
		&$t_{(1,1)}^1$ & $-10 x_1-10 x_2-10 x_3-10 x_4-10 x_5$ \\
		&$t_{(1,1)}^2$ & $x_1 x_2+x_1 x_3+x_2 x_3+x_1 x_4+x_2 x_4+x_3 x_4+x_1 x_5+x_2 x_5+x_3 x_5+x_4 x_5$ \\
		&$t_{(2)}^0$ & $175$ \\
		&$t_{(2)}^1$ & $-21 x_1-21 x_2-21 x_3-21 x_4-21 x_5$ \\
		&$t_{(2)}^2$ & $x_1^2+x_1 x_2+x_2^2+x_1 x_3+x_2 x_3+x_3^2+x_1 x_4+x_2 x_4+x_3 x_4+x_4^2+x_1 x_5+x_2 x_5+x_3 x_5+x_4 x_5+x_5^2$ \\
		5 &$t_{(1,1,1)}^0$ & $-90$ \\
		&$t_{(1,1,1)}^1$ & $25 x_1+25 x_2+25 x_3+25 x_4+25 x_5$ \\
		&$t_{(1,1,1)}^2$ & $-6 x_1 x_2-6 x_1 x_3-6 x_2 x_3-6 x_1 x_4-6 x_2 x_4-6 x_3 x_4-6 x_1 x_5-6 x_2 x_5-6 x_3 x_5-6 x_4 x_5$ \\
		&$t_{(1,1,1)}^3$ & $x_1 x_2 x_3+x_1 x_2 x_4+x_1 x_3 x_4+x_2 x_3 x_4+x_1 x_2 x_5+x_1 x_3 x_5+x_2 x_3 x_5+x_1 x_4 x_5+x_2 x_4 x_5+x_3 x_4 x_5$ \\
		&$t_{(2,1)}^0$  &$-1015$ \\
		&$t_{(2,1)}^1$  &$210 x_1+210 x_2+210 x_3+210 x_4+210 x_5$ \\
		&$t_{(2,1)}^2$  &$-10 x_1^2-31 x_1 x_2-10 x_2^2-31 x_1 x_3-31 x_2 x_3-10 x_3^2-31 x_1 x_4-31 x_2 x_4-31 x_3 x_4-10 x_4^2-31 x_1 x_5-31 x_2 x_5-31 x_3 x_5-31 x_4 x_5-10 x_5^2$ \\
		&$t_{(2,1)}^3$  &$x_1^2 x_2+x_1 x_2^2+x_1^2 x_3+2 x_1 x_2 x_3+x_2^2 x_3+x_1 x_3^2+x_2 x_3^2+x_1^2 x_4+2 x_1 x_2 x_4+x_2^2 x_4+2 x_1 x_3 x_4+2 x_2 x_3 x_4+x_3^2 x_4+x_1 x_4^2+x_2 x_4^2+x_3 x_4^2+x_1^2 x_5+2 x_1 x_2 x_5+x_2^2 x_5+2 x_1 x_3 x_5+2 x_2 x_3 x_5+x_3^2 x_5+2 x_1 x_4 x_5+2 x_2 x_4 x_5+2 x_3 x_4 x_5+x_4^2 x_5+x_1 x_5^2+x_2 x_5^2+x_3 x_5^2+x_4 x_5^2$ \\
		&$t_{(3)}^0$ & $-1960$  \\
		&$t_{(3)}^1$ & $322 x_1+322 x_2+322 x_3+322 x_4+322 x_5$  \\
		&$t_{(3)}^2$ & $-28 x_1^2-28 x_1 x_2-28 x_2^2-28 x_1 x_3-28 x_2 x_3-28 x_3^2-28 x_1 x_4-28 x_2 x_4-28 x_3 x_4-28 x_4^2-28 x_1 x_5-28 x_2 x_5-28 x_3 x_5-28 x_4 x_5-28 x_5^2$  \\
		&$t_{(3)}^3$ & $x_1^3+x_1^2 x_2+x_1 x_2^2+x_2^3+x_1^2 x_3+x_1 x_2 x_3+x_2^2 x_3+x_1 x_3^2+x_2 x_3^2+x_3^3+x_1^2 x_4+x_1 x_2 x_4+x_2^2 x_4+x_1 x_3 x_4+x_2 x_3 x_4+x_3^2 x_4+x_1 x_4^2+x_2 x_4^2+x_3 x_4^2+x_4^3+x_1^2 x_5+x_1 x_2 x_5+x_2^2 x_5+x_1 x_3 x_5+x_2 x_3 x_5+x_3^2 x_5+x_1 x_4 x_5+x_2 x_4 x_5+x_3 x_4 x_5+x_4^2 x_5+x_1 x_5^2+x_2 x_5^2+x_3 x_5^2+x_4 x_5^2+x_5^3$  \\
		& $t_{(1,1,1,1)}^0$ & $31$\\
		& $t_{(1,1,1,1)}^1$ & $-15 x_1-15 x_2-15 x_3-15 x_4-15 x_5$\\
		& $t_{(1,1,1,1)}^2$ & $7 x_1 x_2+7 x_1 x_3+7 x_2 x_3+7 x_1 x_4+7 x_2 x_4+7 x_3 x_4+7 x_1 x_5+7 x_2 x_5+7 x_3 x_5+7 x_4 x_5$\\
		& $t_{(1,1,1,1)}^3$ & $-3 x_1 x_2 x_3-3 x_1 x_2 x_4-3 x_1 x_3 x_4-3 x_2 x_3 x_4-3 x_1 x_2 x_5-3 x_1 x_3 x_5-3 x_2 x_3 x_5-3 x_1 x_4 x_5-3 x_2 x_4 x_5-3 x_3 x_4 x_5$\\
		& $t_{(1,1,1,1)}^4$ & $x_1 x_2 x_3 x_4+x_1 x_2 x_3 x_5+x_1 x_2 x_4 x_5+x_1 x_3 x_4 x_5+x_2 x_3 x_4 x_5$\\
		&$t_{(2,2)}^0$ & $3850$\\
		&$t_{(2,2)}^1$ & $-1050 x_1-1050 x_2-1050 x_3-1050 x_4-1050 x_5$\\
		&$t_{(2,2)}^2$ & $85 x_1^2+225 x_1 x_2+85 x_2^2+225 x_1 x_3+225 x_2 x_3+85 x_3^2+225 x_1 x_4+225 x_2 x_4+225 x_3 x_4+85 x_4^2+225 x_1 x_5+225 x_2 x_5+225 x_3 x_5+225 x_4 x_5+85 x_5^2$\\
		&$t_{(2,2)}^3$ & $-15 x_1^2 x_2-15 x_1 x_2^2-15 x_1^2 x_3-30 x_1 x_2 x_3-15 x_2^2 x_3-15 x_1 x_3^2-15 x_2 x_3^2-15 x_1^2 x_4-30 x_1 x_2 x_4-15 x_2^2 x_4-30 x_1 x_3 x_4-30 x_2 x_3 x_4-15 x_3^2 x_4-15 x_1 x_4^2-15 x_2 x_4^2-15 x_3 x_4^2-15 x_1^2 x_5-30 x_1 x_2 x_5-15 x_2^2 x_5-30 x_1 x_3 x_5-30 x_2 x_3 x_5-15 x_3^2 x_5-30 x_1 x_4 x_5-30 x_2 x_4 x_5-30 x_3 x_4 x_5-15 x_4^2 x_5-15 x_1 x_5^2-15 x_2 x_5^2-15 x_3 x_5^2-15 x_4 x_5^2$\\
		&$t_{(2,2)}^4$ & $x_1^2 x_2^2+x_1^2 x_2 x_3+x_1 x_2^2 x_3+x_1^2 x_3^2+x_1 x_2 x_3^2+x_2^2 x_3^2+x_1^2 x_2 x_4+x_1 x_2^2 x_4+x_1^2 x_3 x_4+2 x_1 x_2 x_3 x_4+x_2^2 x_3 x_4+x_1 x_3^2 x_4+x_2 x_3^2 x_4+x_1^2 x_4^2+x_1 x_2 x_4^2+x_2^2 x_4^2+x_1 x_3 x_4^2+x_2 x_3 x_4^2+x_3^2 x_4^2+x_1^2 x_2 x_5+x_1 x_2^2 x_5+x_1^2 x_3 x_5+2 x_1 x_2 x_3 x_5+x_2^2 x_3 x_5+x_1 x_3^2 x_5+x_2 x_3^2 x_5+x_1^2 x_4 x_5+2 x_1 x_2 x_4 x_5+x_2^2 x_4 x_5+2 x_1 x_3 x_4 x_5+2 x_2 x_3 x_4 x_5+x_3^2 x_4 x_5+x_1 x_4^2 x_5+x_2 x_4^2 x_5+x_3 x_4^2 x_5+x_1^2 x_5^2+x_1 x_2 x_5^2+x_2^2 x_5^2+x_1 x_3 x_5^2+x_2 x_3 x_5^2+x_3^2 x_5^2+x_1 x_4 x_5^2+x_2 x_4 x_5^2+x_3 x_4 x_5^2+x_4^2 x_5^2$\\
		&$t_{(3,1)}^0$ &  $12831$\\
		&$t_{(3,1)}^1$ &  $-3220 x_1-3220 x_2-3220 x_3-3220 x_4-3220 x_5$\\
		&$t_{(3,1)}^2$ &  $280 x_1^2+602 x_1 x_2+280 x_2^2+602 x_1 x_3+602 x_2 x_3+280 x_3^2+602 x_1 x_4+602 x_2 x_4+602 x_3 x_4+280 x_4^2+602 x_1 x_5+602 x_2 x_5+602 x_3 x_5+602 x_4 x_5+280 x_5^2$\\
		&$t_{(3,1)}^3$ &  $-10 x_1^3-38 x_1^2 x_2-38 x_1 x_2^2-10 x_2^3-38 x_1^2 x_3-66 x_1 x_2 x_3-38 x_2^2 x_3-38 x_1 x_3^2-38 x_2 x_3^2-10 x_3^3-38 x_1^2 x_4-66 x_1 x_2 x_4-38 x_2^2 x_4-66 x_1 x_3 x_4-66 x_2 x_3 x_4-38 x_3^2 x_4-38 x_1 x_4^2-38 x_2 x_4^2-38 x_3 x_4^2-10 x_4^3-38 x_1^2 x_5-66 x_1 x_2 x_5-38 x_2^2 x_5-66 x_1 x_3 x_5-66 x_2 x_3 x_5-38 x_3^2 x_5-66 x_1 x_4 x_5-66 x_2 x_4 x_5-66 x_3 x_4 x_5-38 x_4^2 x_5-38 x_1 x_5^2-38 x_2 x_5^2-38 x_3 x_5^2-38 x_4 x_5^2-10 x_5^3$\\
		&$t_{(3,1)}^4$ &  $x_1^3 x_2+x_1^2 x_2^2+x_1 x_2^3+x_1^3 x_3+2 x_1^2 x_2 x_3+2 x_1 x_2^2 x_3+x_2^3 x_3+x_1^2 x_3^2+2 x_1 x_2 x_3^2+x_2^2 x_3^2+x_1 x_3^3+x_2 x_3^3+x_1^3 x_4+2 x_1^2 x_2 x_4+2 x_1 x_2^2 x_4+x_2^3 x_4+2 x_1^2 x_3 x_4+3 x_1 x_2 x_3 x_4+2 x_2^2 x_3 x_4+2 x_1 x_3^2 x_4+2 x_2 x_3^2 x_4+x_3^3 x_4+x_1^2 x_4^2+2 x_1 x_2 x_4^2+x_2^2 x_4^2+2 x_1 x_3 x_4^2+2 x_2 x_3 x_4^2+x_3^2 x_4^2+x_1 x_4^3+x_2 x_4^3+x_3 x_4^3+x_1^3 x_5+2 x_1^2 x_2 x_5+2 x_1 x_2^2 x_5+x_2^3 x_5+2 x_1^2 x_3 x_5+3 x_1 x_2 x_3 x_5+2 x_2^2 x_3 x_5+2 x_1 x_3^2 x_5+2 x_2 x_3^2 x_5+x_3^3 x_5+2 x_1^2 x_4 x_5+3 x_1 x_2 x_4 x_5+2 x_2^2 x_4 x_5+3 x_1 x_3 x_4 x_5+3 x_2 x_3 x_4 x_5+2 x_3^2 x_4 x_5+2 x_1 x_4^2 x_5+2 x_2 x_4^2 x_5+2 x_3 x_4^2 x_5+x_4^3 x_5+x_1^2 x_5^2+2 x_1 x_2 x_5^2+x_2^2 x_5^2+2 x_1 x_3 x_5^2+2 x_2 x_3 x_5^2+x_3^2 x_5^2+2 x_1 x_4 x_5^2+2 x_2 x_4 x_5^2+2 x_3 x_4 x_5^2+x_4^2 x_5^2+x_1 x_5^3+x_2 x_5^3+x_3 x_5^3+x_4 x_5^3$\\
		&$t_{(4)}^0$ & $22449$ \\
		&$t_{(4)}^1$ & $-4536 x_1-4536 x_2-4536 x_3-4536 x_4-4536 x_5$ \\
		&$t_{(4)}^2$ & $546 x_1^2+546 x_1 x_2+546 x_2^2+546 x_1 x_3+546 x_2 x_3+546 x_3^2+546 x_1 x_4+546 x_2 x_4+546 x_3 x_4+546 x_4^2+546 x_1 x_5+546 x_2 x_5+546 x_3 x_5+546 x_4 x_5+546 x_5^2$ \\
		&$t_{(4)}^3$ & $-36 x_1^3-36 x_1^2 x_2-36 x_1 x_2^2-36 x_2^3-36 x_1^2 x_3-36 x_1 x_2 x_3-36 x_2^2 x_3-36 x_1 x_3^2-36 x_2 x_3^2-36 x_3^3-36 x_1^2 x_4-36 x_1 x_2 x_4-36 x_2^2 x_4-36 x_1 x_3 x_4-36 x_2 x_3 x_4-36 x_3^2 x_4-36 x_1 x_4^2-36 x_2 x_4^2-36 x_3 x_4^2-36 x_4^3-36 x_1^2 x_5-36 x_1 x_2 x_5-36 x_2^2 x_5-36 x_1 x_3 x_5-36 x_2 x_3 x_5-36 x_3^2 x_5-36 x_1 x_4 x_5-36 x_2 x_4 x_5-36 x_3 x_4 x_5-36 x_4^2 x_5-36 x_1 x_5^2-36 x_2 x_5^2-36 x_3 x_5^2-36 x_4 x_5^2-36 x_5^3$ \\
		&$t_{(4)}^4$ & $x_1^4+x_1^3 x_2+x_1^2 x_2^2+x_1 x_2^3+x_2^4+x_1^3 x_3+x_1^2 x_2 x_3+x_1 x_2^2 x_3+x_2^3 x_3+x_1^2 x_3^2+x_1 x_2 x_3^2+x_2^2 x_3^2+x_1 x_3^3+x_2 x_3^3+x_3^4+x_1^3 x_4+x_1^2 x_2 x_4+x_1 x_2^2 x_4+x_2^3 x_4+x_1^2 x_3 x_4+x_1 x_2 x_3 x_4+x_2^2 x_3 x_4+x_1 x_3^2 x_4+x_2 x_3^2 x_4+x_3^3 x_4+x_1^2 x_4^2+x_1 x_2 x_4^2+x_2^2 x_4^2+x_1 x_3 x_4^2+x_2 x_3 x_4^2+x_3^2 x_4^2+x_1 x_4^3+x_2 x_4^3+x_3 x_4^3+x_4^4+x_1^3 x_5+x_1^2 x_2 x_5+x_1 x_2^2 x_5+x_2^3 x_5+x_1^2 x_3 x_5+x_1 x_2 x_3 x_5+x_2^2 x_3 x_5+x_1 x_3^2 x_5+x_2 x_3^2 x_5+x_3^3 x_5+x_1^2 x_4 x_5+x_1 x_2 x_4 x_5+x_2^2 x_4 x_5+x_1 x_3 x_4 x_5+x_2 x_3 x_4 x_5+x_3^2 x_4 x_5+x_1 x_4^2 x_5+x_2 x_4^2 x_5+x_3 x_4^2 x_5+x_4^3 x_5+x_1^2 x_5^2+x_1 x_2 x_5^2+x_2^2 x_5^2+x_1 x_3 x_5^2+x_2 x_3 x_5^2+x_3^2 x_5^2+x_1 x_4 x_5^2+x_2 x_4 x_5^2+x_3 x_4 x_5^2+x_4^2 x_5^2+x_1 x_5^3+x_2 x_5^3+x_3 x_5^3+x_4 x_5^3+x_5^4$ \\
		\hline
\end{longtable}
%	\end{tabular}
\end{center}

%%%%%%%%%%%%%%%%%%%%%%%%%%%%%%%% Acknowledgement%%%%%%%
{\bf Acknowledgements} We are indebted to D. Trusin, M. Movshev, M. Mulase and F. Plaza-Martin for useful comments. Our special thanks to B. Osserman and V. Vologodsky for help with  generalization of our results to integral curves.

\bibliographystyle{amsplain}
%\bibliography{sampartb}

\end{document}